\documentclass[12pt]{amsart}

\makeatletter

\renewcommand{\subsection}{\@startsection{subsection}{2}%
  \z@{.5\linespacing\@plus.7\linespacing}{-.5em}%
  {\normalfont\bfseries\S\,}}

\makeatother

\newcommand{\qsp}{\gQ^\sep_k}
\newcommand{\csp}{\gC^\sep_k}

\usepackage{amsmath}
\usepackage{amssymb}
\usepackage{enumerate}
\usepackage{euscript}
\usepackage{amscd}
\usepackage[all]{xy}
\newtheorem{thm}{Theorem}[section]
\newtheorem{lem}[thm]{Lemma}
\newtheorem{prop}[thm]{Proposition}

\theoremstyle{definition}
\newtheorem{defn}[thm]{Definition}
\newtheorem{rem}[thm]{Remark}

\newtheorem{asmp}[thm]{Assumption}

\numberwithin{equation}{section}
\theoremstyle{remark}

\newcommand{\bba}{{\mathbb A}}

\newcommand{\bbc}{{\mathbb C}}

\newcommand{\bbq}{{\mathbb Q}}
\newcommand{\bbr}{{\mathbb R}}

\newcommand{\bbz}{{\mathbb Z}}

\newcommand{\lam}{{\lambda}}
\newcommand{\Lam}{{\Lambda}}
\newcommand{\lamb}{\underline{\lambda}}
\newcommand{\mub}{\underline{\mu}}

\newcommand{\gc}{{\mathfrak c}}
\newcommand{\gC}{{\mathfrak C}}

\newcommand{\gM}{{\mathfrak M}}

\newcommand{\gQ}{{\mathfrak Q}}

\newcommand{\gS}{{\mathfrak S}}

\newcommand{\gV}{{\mathfrak V}}

\newcommand{\cK}{{\mathcal K}}

\newcommand{\cR}{{\mathcal R}}
\newcommand{\co}{{\mathcal O}}

\newcommand{\gl}{{\operatorname{GL}}}

\newcommand{\Res}{{\operatorname{Res}}}

\newcommand{\op}{{\operatorname{op}}}

\newcommand{\sst}{{\operatorname{ss}}}

\newcommand{\sep}{{\operatorname{sep}}}

\newcommand{\psf}{Poisson summation formula}

\newcommand{\res}{\operatorname{Res}}

\newcommand{\A}{\bba}
\newcommand{\Z}{\bbz}
\newcommand{\Q}{\bbq}
\newcommand{\R}{\bbr}
\newcommand{\C}{\bbc}

\newcommand{\ma}{\bba^{\times}}

\newcommand{\mk}{k^{\times}}
\newcommand{\ag}{G_{\bba}}
\newcommand{\rg}{G_{k}}

\newcommand{\av}{V_{\bba}}

\newcommand{\md}{d^{\times}\!}
\newcommand{\gMf}{\gM_{\text f}}
\newcommand{\gMc}{\gM_\C}
\newcommand{\gMr}{\gM_\R}

\newcommand{\ti}{\widetilde}

\newcommand{\ct}{\frac 1 {2\pi\sqrt{-1}}}

\newcommand{\sE}{\mathcal {E}}
\newcommand{\sS}{\mathcal {S}}

\newcommand{\ac}[1]{\langle{#1}\rangle}% additive character
\newcommand{\nr}{\EuScript N}
\newcommand{\trc}{\EuScript T}
\newcommand{\inv}{{\operatorname{Inv}}}
\newcommand{\rnk}{{\operatorname{rank}}}
\renewcommand{\cK}{K}

\begin{document}

\title[pair of simple algebras]
{On the zeta functions of prehomogeneous vector spaces
for pair of simple algebras}
%$\text{\Large for the space of pairs of quaternions}$}
\author[Takashi Taniguchi]{Takashi Taniguchi}
\keywords{prehomogeneous vector space, global zeta function, simple algebra,
density theorem}
\address{Graduate School of Mathematical Sciences\\ University of Tokyo\\
3--8--1 Komaba Megoro-ku\\ Tokyo 153-0041\\ JAPAN}
\email{tani@ms.u-tokyo.ac.jp}
\date{\today}

\begin{abstract}
In this paper, we consider the
prehomogeneous vector space for pair of simple algebras
which are $k$-forms of the $D_4$ type and the $E_6$ type.
We mainly study the non-split cases.
The main purpose of this paper is to
determine the principal parts of the global zeta functions
associated with these spaces in the non-split cases.
We also give a description of the sets of rational orbits of these spaces,
which suggests the expected density theorems arising from the properties of 
these zeta functions.
\end{abstract}
%%%%%%%%%%	TEXT START		%%%%%%%%%%

\maketitle
%\setcounter{tocdepth}{1}
%\tableofcontents	

%%%%%%%%%% Section	%%%%%%%%%%
\section{Introduction}

Let $k$ be a field and $D$ a simple algebra
of dimension $4$ or $9$ over $k$.
We denote by $D^{\rm op}$ the opposite algebra of $D$.
In this paper, we consider the prehomogeneous vector space
$(G,V)=(G,\rho,V)$ where
\begin{equation}\label{ourpv}
G=D^\times\times(D^\op)^\times \times\gl(2),\qquad V=D\otimes k^2,
\end{equation}
and
\begin{equation*}
\rho(g)(a\otimes v)=(g_{11}ag_{12})\otimes (g_2v)
\quad\text{for}\quad
g=(g_{11},g_{12},g_2)\in G, a\in D, v\in k^2.
\end{equation*}
We call the prehomogeneous vector space $D_4$ type and $E_6$ type
if the dimension of $D$ is $4$ and $9$, respectively.
This representation is a $k$-form of
\begin{equation*}
G'=\gl(n)\times \gl(n)\times \gl(2),\qquad V'=k^n\otimes k^n\otimes k^2
\end{equation*}
for $n=2$ and $n=3$ if the dimension of $D$ is $4$ and $9$, respectively.
If $D$ is split, then $(G,V)$ is equivariant to $(G',V')$ over $k$.
In this paper, we give a certain description of $k$-rational orbits
and determine the structure of the stabilizers for semi-stable points.
Also we determine the principal parts of the global zeta function
for the non-split cases of $(G,V)$ over an algebraic number field $k$.

We recall the definition of prehomogeneous vector spaces.
Here we give a definition of a certain restricted class
instead of giving general one for simplicity.

\begin{defn}\label{pv}
An irreducible representation of a connected reductive group
$(G,V)$ over $k$ is called a {\em prehomogeneous vector space} if
\begin{enumerate}[{\rm (1)}]
\item there exists a Zariski open $G$-orbit in $V$ and
\item there exists a non-constant polynomial $P\in k[V]$
and a rational character $\chi$ of $G$ such that
$P(gx)=\chi(g)P(x)$ for all $g\in G$ and $x\in V$.
\end{enumerate}
\end{defn}

Irreducible prehomogeneous vector spaces over an arbitrary
characteristic $0$ algebraically closed field
was classified by Sato and Kimura in \cite{saki}.
Sato and Shintani \cite{sash} defined
global zeta functions for prehomogeneous vector spaces
if $(G,V)$ is defined over a number field.

The information of the principal part at the rightmost pole
of the global zeta function for a prehomogeneous vector space
together with an appropriate local theory
yields interesting density theorems.
For example, 
using Shintani's result \cite{shintania}
for the space of binary cubic forms,
Datskovsky and Wright \cite{dawra,dawrb}
gave the zeta function theoretic proof
of the Davenport and Heilbronn \cite{dahea,daheb} density theorem
\begin{equation*}
\sum_{\underset{|\Delta_F|\leq x}{[F:\Q]=3}}1\sim \frac{x}{\zeta(3)}
\qquad (x\to\infty),
\end{equation*}
where $F$ runs through all the cubic fields
with the absolute value of its
discriminant $|\Delta_F|$ is not bigger than $x$.
Also recent works concerning to the space of pairs of binary Hermitian forms
\cite{kable-yukie-pbh-I,kable-yukie-pbh-II,kable-yukie-pbh-III}
by Kable and Yukie
gave certain new density theorems
together with Yukie's global theory \cite{yukieh}.
For the statement of the density theorem, see the introduction of
\cite{kable-yukie-pbh-I}.
Note that this case is another $k$-form of the $D_4$ case.
These $k$-forms are listed in H. Saito's classification \cite{hsaitoa}.

We return to our prehomogeneous vector space \eqref{ourpv}.
The following theorem is a main result of this paper.
%In our case, we obtain the following theorem
%information of the rightmost 
%pole of the global zeta fucntion 
%as a corollary to Theorem \ref{ppf}.
\begin{thm}\label{rmp}
Let $D$ be a non-split simple algebra of dimension $m=4$ or $9$.
The global zeta function $Z(\Phi,s)$ 
associated with the prehomogeneous vector space \eqref{ourpv}
can be continued meromorphically
to the region $\Re(s)>2m-2$ only with a possible simple pole
at $s=2m$ with the residue
$\tau(G_1)\gV_2\int_{V_\A}\Phi(x)dx$.
\end{thm}
%The global zeta function $Z(\Phi,s)$,
The constants $\tau(G_1),\gV_2$
and the measure $dx$ on $V_\A$ are defined in Section \ref{gzf}.
All the other poles are also
described by means of certain distributions in Theorem \ref{ppf},
but the above theorem is enough to get the density theorems.
On the other hand, the expected density theorems require
not only the standard tauberian theorem,
but also the appropriate local theory and
what is called ``filtering process'',
which will be studied in a forthcoming paper.
For a general transition process from 
the tauberian theorem for global zeta function
of prehomogeneous vector spaces to density theorems, 
see \cite{yukiec}.

The expected density theorems from our cases will be discussed
in Remark \ref{exp} using a result in Section \ref{rod},
but we also give a brief summary here.
We assume $k=\Q$ for simplicity.
For a finite extension $F$ of $\Q$, $h_F$ and $R_F$ denote
the class number and the regulator of $F$, respectively.

For a prehomogeneous vector space $(G,V)$ in Definition \ref{pv},
set $V^\sst_\Q=\{x\in V_\Q\mid P(x)\not=0\}$ where $P(x)$ is as in (2)
and $\ti T=\ker(G\rightarrow \gl(V))$.
For $x\in V_\Q^\sst$, $G_x$ denotes the stabilizer of $x$
and $G_x^\circ$ denotes its identity component.
Roughly speaking, the global zeta function is a counting function
for the unnormalized Tamagawa numbers of $G_x^\circ/\ti T$ of points in
$x\in G_\Q\backslash V_\Q^\sst$ (see \cite{yukiec}, for example).
Hence, by Proposition \ref{rod-d4} and \ref{rod-e6},
the behavior of the zeta function at the rightmost pole
for the $D_4$ case and the $E_6$ case
(both split and non-split)
yields the asymptotic behavior of the following function
\begin{equation*}
\sum_{\underset{\text{$F$ is embeddable into $D$}}{\underset{|\Delta_F|\leq x}{[F:\Q]=2}}}\hspace*{-7mm}h_F^2R_F^2
\qquad\text{and}\quad
\sum_{\underset{\text{$F$ is embeddable into $D$}}{\underset{|\Delta_F|\leq x}{[F:\Q]=3}}}\hspace*{-7mm}h_FR_F
\end{equation*}
as $x\to\infty$, respectively.
If the method of filtering process in \cite{dawra, dawrb},
\cite{dats},
\cite{kable-yukie-pbh-I,kable-yukie-pbh-II,kable-yukie-pbh-III}
is applicable to compute the above density,
one can also compute a similar density with 
finitely many local conditions
by the same method.
%Moreover, if one succeeds to
%establish the above density theorems,
%he can also obtain the density with $F$ runs through
%given local behavior at a fixed finite number of places.
It is also easy to show that for each $D$, the condition $F\subset D$
can be determined by finitely many local conditions of $F$.
Hence, the expected density theorems from the non-split cases
seems to be quite similar to that of from the split cases.

One advantage of non-split cases is
that the global theory becomes much easier.
%As is occured in \cite{yukiec}, 
The analysis of the global zeta function becomes much complicated
as the split rank of the group growth,
and we could not yet success to establish the global theory
for the split $D_4$ and $E_6$ cases.
Especially the rank of the group is $5$ for the split $E_6$ case,
and the complexity of computing the principal part of the zeta function
seems to be formidable.

For the rest of this section, we will give the contents of this paper
and the notations used in this paper.
More specific notations will be introduced in each section.
In Section \ref{psa}, we will define the space of pair of simple algebras
and summarize its basic properties.
Before starting the global theory in Section \ref{gzf},
we will give a certain description
of rational orbits in Section \ref{rod}.
The split cases are treated in \cite[\S3]{wryu}, and
this is a slight generalization of those cases.
We prove that the set of rational orbits corresponds
one to one to a certain set
of quadratic extensions and cubic extensions
for the $D_4$ case and the $E_6$ case, respectively.
Also we determine the structure of the stabilizers for semi-stable points. 
The expected density theorems from our cases will be discussed
in Remark \ref{exp}.

In Section \ref{gzf}, we study the global theory for the non-split cases.
In \S \ref{ss:pre}, we introduce notations used in this section and
review some basic facts on adelic analysis.
In \S \ref{ss:gzf}, we define a global zeta function.
Also we will give an estimate of an incomplete theta series.
Although H. Saito \cite{hsaitoc} proved the convergence of the 
all global zeta functions
associated with prehomogeneous vector spaces, 
we need the estimate in order to use Shintani's lemma.
In \S \ref{ss:pp}, we divide the global zeta function into
the ``entire part'' and the ``principal part'' 
%describe the principal part of the zeta function
by using the Poisson summation formula.
We study the ``principal part'' in later subsections.

In \S \ref{ss:str}, we introduce a stratification of unstable points.
To separate the contribution from unstable strata, we use Shintani's lemma.
In \S \ref{ss:ses}, we review Shintani's lemma and apply it to our cases.
Since $D^\times$ is of rank $0$, 
the smoothed Eisenstein series for our case is
essentially the same as that of $\gl(2)$.
In \S \ref{ss:zsa}, we review some analytic properties of 
the zeta function associated with (single) simple algebra,
because this zeta function appears in the induction process.
In \S \ref{ss:cus}, we compute contributions from unstable points.
By putting the results together we have obtained
in \S \ref{ss:str}--\S \ref{ss:cus},
we determine the principal part of the global zeta function
in \S \ref{ss:ppf}.

The standard symbols $\Q$, $\R$, $\C$ and $\Z$ will denote respectively
the rational, real and complex numbers and the rational integers.
If $R$ is any ring then $R^{\times}$ is the set
of invertible elements of $R$ and if $V$ is a variety defined over
$R$ then $V_R$ denotes its $R$-points.
In Section \ref{psa} and \ref{rod}, $k$ denotes arbitrary field.
In Section \ref{gzf}, $k$ denotes an algebraic number field.

%\vspace{1zw}

\noindent
{\bf Acknowledgments.}
This work is the author's doctor's thesis at University of Tokyo.
The author would like to express his sincere gratitude
to his advisor T. Terasoma
for the constant support, encouragement and many helpful suggestions.
The author would like to heartily thank Professor A. Yukie
for the encouragement and many invaluable inspiring suggestions.
The author is deeply grateful to Professor T. Oda
for the support and consultations.
He also would like to thank his colleague T. Ito
with the many useful discussions.
This work is partially supported by
21st Century (University of Tokyo) COE program,
of Ministry of Education, Culture, Sports, Science and Technology.

%\newpage
\section{The space of pair of simple algebras}\label{psa}
In this section, we define the representation
of the space of pair of simple algebras,
and discuss its basic properties.

Let $k$ be an arbitrary field and $D$ a simple algebra over $k$
of dimension $m=n^2,n\geq1$.
Let $\trc$ and $\nr$ be the reduced trace and reduced norm, respectively.
We denote by $D^\op$ the opposite algebra of $D$.
We introduce a group $G_1$ and its representation space $W$ as follows.
Let
\begin{equation*}
G_1=D^\times\times(D^\op)^\times.
\end{equation*}
That is, $G_1$ is equal to $D^\times\times D^\times$ set theoretically
and the multiplication law is given by
$(g_{11},g_{12})(h_{11},h_{12})=(g_{11}h_{11},h_{12}g_{12})$.
We regard $G_1$ as an algebraic group over $k$.
The simple algebra $D$ can be considered
as a vector space over $k$.
When we regard $D$ as a vector space over $k$,
we denote this space as $W$.
We define the action of $G_1$ on $W$ as follows:
\begin{equation*}
(g_1,w)\longmapsto g_{11}w g_{12},
	\qquad g_1=(g_{11},g_{12})\in G_1,
			w\in W.
\end{equation*}
This defines a representation $W$ of $G_1$.
Clearly, $(G_1,W)$ is a prehomogeneous vector space.
We discuss the properties of the zeta function
associated with this space in \S\ref{ss:zsa},
which will be used in the analysis of the zeta function
associated with the space of pair of simple algebras.

Let $G_2=\gl(2)$ and $k^2$ the standard representation of $G_2$.
The group $G=G_1\times G_2$ acts naturally on
$V=W\otimes k^2$.
This is a $k$-form of
$(\gl(n)\times\gl(n)\times\gl(2),k^n\otimes k^n\otimes k^2)$,
and it is proved in \cite{saki} that 
this is a prehomogeneous vector space 
if and only if $n=2$ or $n=3$.
Since we are interested in prehomogeneous vector space,
{\em we consider the case $n=2,3$ for the rest of this paper}.
That is, $D$ is a simple algebra of dimension $4$ or $9$.
We call these representation $D_4$ type and $E_6$ type
for $n=2$ and $n=3$, respectively following \cite{wryu}.

We describe the action more explicitly.
Throughout of this paper,
we express elements of $V\cong W\oplus W$ as $x=(x_1,x_2)$.
We identify $x=(x_1,x_2)\in V$ with $x(v)=v_1x_1+v_2x_2$
which is an element of simple algebra with entries
in linear forms in two variables $v=(v_1,v_2)$.
Then the action of $g=(g_{11},g_{12},g_2)\in G$ on $x\in V$ is defined by
\begin{equation*}
(gx)(v)=g_{11}x(vg_2)g_{12}.
\end{equation*}
We put $F_x(v)=\nr(x(v))$.
This is a binary quadratic form (resp.\ cubic form)
in variables $v=(v_1,v_2)$ if $n=2$ (resp.\ $n=3$),
and the discriminant $P(x)\ (x\in V)$ is a polynomial in $V$.
The polynomial $P(x)$ is characterized by 
\begin{equation*}
P(x)=\prod_{i<j}(\alpha_i\beta_j-\alpha_j\beta_i)^2
\quad
\text {for}
\quad
F_x(v)=\prod_{1\leq i\leq n}(\alpha_iv_1-\beta_iv_2), \ \ x\in V_{\bar k}.
\end{equation*}
Let $\chi_i\, (i=1,2)$ be the character of $G_i$ defined by
\begin{equation*}
\chi_1(g_1)=\nr(g_{11})\nr(g_{12}),\quad \chi_2(g_2)=\det g_2,
\end{equation*}
respectively.
We define $\chi(g)=\chi_1(g_1)^2\chi_2(g_2)^2$ for $n=2$ and
$\chi(g)=\chi_1(g_1)^4\chi_2(g_2)^6$ for $n=3$.
Then one can easily see that
\begin{equation*}
P(gx)=\chi(g)P(x)
\end{equation*}
and hence $P(x)$ is a relative invariant polynomial
with respect to the character $\chi$.
Let $S=\{x\in V\mid P(x)=0\}$ and $V^\sst=\{x\in V\mid P(x)\not=0\}$
and call them the set of unstable points and semi-stable points, respectively.
That is, $x\in V$ is semi-stable if and only if
$F_x(v)$ does not have a multiple root in ${\mathbb P}^1=\{(v_1:v_2)\}$.

%\newpage
\section{Rational orbit decomposition}\label{rod}
\subsection{Rational orbit decomposition}
%As we mentioned in the introduction, .........
In this section, we will interpret the rational orbit space
$G_k\backslash V_k^\sst$ and determine the structure of the
stabilizers for semi-stable points.
The split cases are treated in \cite[\S3]{wryu},
and here is a slight generalization of that.
For the expected density theorems, see Remark \ref{exp}.
For $x\in V_k^\sst$, let $G_x$ be the stabilizer of $x$
and $G_x^\circ$ its identity component.

Let $\qsp$ (resp.\ $\csp$) be set of isomorphism classes
of separable commutative $k$-algebras of dimension $2$ (resp.\ $3$).
For example, $\qsp$ can be regarded as
the disjoint union of $\{k^2\}$ and the set of
separable quadratic extensions of $k$.
\begin{defn}
For $x\in V_k^\sst$, we define
\begin{align*}
Z_x		&=\mathrm{Proj}\, k[v_1,v_2]/(F_x(v)),\\
\ti k(x)	&=\Gamma(Z_x,\co_{Z_x}).
\end{align*}
Also we define $k(x)$ to be the splitting field of $F_x(v)$.
\end{defn}
Note that $\ti k(x)$ may not be a field.
Since $V_k^\sst$ is the set of $x$ such that
$F_x$ does not have a multiple root,
$Z_x$ is a reduced scheme over $k$ and
$\ti k(x)$ is an element of
$\qsp$ (resp.\ $\csp$) for $n=2$ (resp.\ $n=3$).
%For $x\in V_k^\sst$, .............
Since
\begin{equation*}
F_{gx}(v)=\chi_1(g_1)F_x(vg_2),
\end{equation*}
the isomorphism classes of
$Z_x, \ti k(x)$ and $k(x)$ depend only on the $G_k$-orbit of $x$.

We let
\begin{equation*}
\ti\alpha_V\colon
	G_k\backslash V_k^\sst	\longrightarrow	\qsp\ \ (\text{resp.}\, \csp)
\quad
	x\longmapsto \ti k(x)
\end{equation*}
for the $D_4$ case (resp.\ the $E_6$ case).
We first determine the image of $\ti \alpha_V$.

\begin{defn}
\begin{enumerate}[{\rm (1)}]
\item For $n=2$, 
we denote by $\qsp(D)$ the subset of $\qsp$
consisting of algebras which have an embedding to $D_k$.
\item For $n=3$, 
we denote by $\csp(D)$ the subset of $\csp$
consisting of algebras which have an embedding to $D_k$.
\end{enumerate}
\end{defn}

\begin{lem}\label{image}
\begin{enumerate}[{\rm (1)}]
\item Let $(G,V)$ be the $D_4$ case.
The image of the map $\ti\alpha_V$ is $\qsp(D)$.
\item Let $(G,V)$ be the $E_6$ case.
The image of the map $\ti\alpha_V$ is $\csp(D)$.
\item Moreover, any orbit $G_kx\subset V_k^\sst$
contains an element of the form $y=(1,y_2)$.
\end{enumerate}
\end{lem}
\begin{proof}
Here we consider the $E_6$ case.
The $D_4$ case can be treated similarly.
First note that for $x=(1,w)\in V_k^\sst$,
$F_x(v_1,1)=\nr(v_1+w)$ is the characteristic polynomial of $-w\in D_k$
that does not have a multiple root,
and hence the algebra
$\ti k(x)$ is isomorphic to the subalgebra $k[w]\subset D_k$
generated by $w$ over $k$ in $D_k$.

Let $L\in \csp(D)$. We regard $L$ as a subalgebra of $D_k$ and
take an element $u\in D_k$ so that $L=k[u]$.
Let $x=(1,-u)\in V_k$.
Then since $\dim_kL=3=\deg F_x(v_1,1)$,
the characteristic polynomial $F_x(v_1,1)$ of $u$
is also the minimum polynomial of $u$.
Hence $F_x(v_1,1)$ does not have a multiple root since $u$ is separable.
This shows that $x=(1,-u)\in V_k^\sst$ and now
by the remark above we have $\ti\alpha_V(G_kx)=L$.
This proves that the image of $\ti\alpha_V$ contains $\csp(D)$.

Since (3) implies the opposite inclusion, we consider (3).
If $D_k$ is non-split, this is obvious because $D_k$ is a division algebra.
We consider the split cases.
Since the argument is similar, we consider the $E_6$ case here.
In this case, $D={\mathrm M}(3)$ be the algebra of $3\times 3$ matrices.
For $a\in D$, let $\rnk(a)$ denote the rank of
the matrix $a$.

Let $x=(x_1,x_2)\in V_k^\sst$.
If either the rank of $x_1$ or $x_2$ is equal to $3$, 
the element is invertible and hence, there exists a $g\in G_{k}$
such that $gx=(1,*)$.
Also if both the rank of $x_1$ and $x_2$ are less than or equal to $1$,
we have $F_x(v)=\det(x_1v_1+x_2v_2)=0$ which contradicts to $x\in V_k^\sst$.
Hence, by changing if necessary, we assume that $\rnk(x_1)=2$.
Then there exists a $g_1\in G_{1k}$ such that $x'=g_1x=(e,y)$, where
\begin{equation*}
e= \begin{pmatrix} 1 & 0 & 0\\ 0 & 1 & 0\\ 0& 0& 0\end{pmatrix},\quad
y= \begin{pmatrix} y_{11} & y_{12} & y_{13}\\ y_{21} & y_{22} & y_{23}\\
 y_{31}& y_{32}& y_{33}\end{pmatrix}.
\end{equation*}
If $y_{33}=0$, then it is easy to see that $F_{x'}(v)$ has a multiple root,
and so $y_{33}\not=0$.
Hence, again we can take an element $g_1'\in G_{1k}$ so that
\begin{equation*}
g_1'(e,y)=(e,z),\quad
z= \begin{pmatrix} z_{11} & z_{12} & 0\\ z_{21} & z_{22} & 0\\
 0& 0& 1\end{pmatrix}.
\end{equation*}
Now it is easy to see that there exist $\alpha,\beta\in k$ so that
$\rnk(\alpha e+\beta z)=3$, hence we have (3).
\end{proof}

We later show that the map $\ti\alpha_V$ is in fact injective.
Next we consider the structure of the stabilizers for semi-stable points.
Note that for $x\in V_k^\sst$,
\begin{equation*}
\dim G_x^\circ=\dim G_x=\dim G-\dim V=4.
\end{equation*}
\begin{lem}
Let $x\in V_k^\sst$.
\begin{enumerate}[{\rm (1)}]
\item Let $(G,V)$ be the $D_4$ case.
 We have $G_x^\circ\cong(\gl(1)_{\ti k(x)})^2$ as a group over $k$.
\item Let $(G,V)$ be the $E_6$ case.
 We have $G_x^\circ\cong\gl(1)_{\ti k(x)}\times\gl(1)_k$ as a group over $k$.
\end{enumerate}
\end{lem}

\begin{proof}
By Lemma \ref{image} (3),
any $G_k$-orbit in $V_k^\sst$ contains an element of the form
$x=(1,w)$ with $w\in D_k\setminus k$.
Hence it is enough to show the lemma for these elements.
We identify $\ti k(x)$ with $k[w]\subset D_k$.

In order to prove an isomorphism
between two algebraic groups $G_1,G_2$ over $k$,
it is enough to construct isomorphisms
between the sets $G_{1R},G_{2R}$ of $R$-rational points
of $G_1,G_2$ for all commutative $k$-algebras $R$
which satisfy the usual functorial property.
For this, the reader should see \cite[p.17]{mum}.

We first consider (1). For the $D_4$ case,
$\ti k(x)=k[w]$ is a separable $k$-algebra of dimension $2$.
Let $R$ be any commutative $k$-algebra. We put $\ti R(x)= \ti k(x)\otimes R$.
Note that $\ti R(x)=R[w]$ is a subalgebra of
$D_R=D\otimes R$ and is commutative.
Since $\{1,w\}$ is a $k$-basis of $\ti k(x)$,
this is also an $R$-basis of $\ti R(x)$.
Let $s,t\in\ti R(x)^\times$.
Then $\{st,stw\}$ is also an $R$-basis of $\ti R(x)$,
and so there exists a unique element $g=g_{st}\in\gl(2)_R$ such that
$g\,\,^{t}(st,stw)={}^t(1,w)$.
Hence $(s,t)\mapsto (s,t,g_{st})$ gives an injective homomorphism
from $(\ti R(x)^\times)^2$ to $G_{xR}$.

This shows that there exists an injective homomorphism
\begin{equation*}
(\gl(1)_{\ti k(x)})^2\longrightarrow G_x.
\end{equation*}
Since $(\gl(1)_{\ti k(x)})^2$ is a connected algebraic group of dimension $4$,
we have $(\gl(1)_{\ti k(x)})^2\cong G_x^\circ$.

Next we consider (2).
Again we let $R$ be any algebra and
put $\ti R(x)=\ti k(x)\otimes R$.
Then we have a injective homomorphism
from $\ti R(x)^\times\times R^\times$ to $G_{xR}$
by sending $(s,t)$ to $(s,s^{-1}t^{-1},t)$.
This shows that there exists an injective homomorphism
\begin{equation*}
\gl(1)_{\ti k(x)}\times \gl(1)_k\longrightarrow G_x.
\end{equation*}
Since $\gl(1)_{\ti k(x)}\times \gl(1)_k$ is a
connected algebraic group of dimension $4$,
we have $\gl(1)_{\ti k(x)}\times\gl(1)_k\cong G_x^\circ$.
\end{proof}

Finally, we show the injectivity of $\ti\alpha_V$.
\begin{lem}
In both cases, the map $\ti\alpha_V$ is injective.
\end{lem}
\begin{proof}
%.......
%The author hopes that if we suitably rewrite the argument of
%\cite{wryu}, \cite{kayu} in our settings,
%we can show this lemma from the fact that
%each Galois cohomology of $G_x^\circ$ vanishes,
%but has not established yet.
%On the other side, if we restrict ourselves to consider
%only for non-split cases,
%there is an alternative elementary proof due to the fact that
%any element of $\qsp(D)$ or $\csp(D)$ is a field.
%Since the split case is already proven in \cite{wryu}, 
%we put the proof at present.......

Since the split case is already proven in \cite{wryu}, 
we only consider the non-split cases here.
Let $x,y\in V_k^\sst$ satisfy $\ti k(x)\cong \ti k(y)$.
By Lemma \ref{image} (3), we may assume $x=(1,u_1),y=(1,u_2)$.
Then $k[u_1]$ and $k[u_2]$ are isomorphic subfields of $D_k$.
By the Skolem-Noether theorem \cite[CHAPITRE 8 \S10]{bralgebre},
there exists an element $\theta\in D_k^\times$ such that
\begin{equation*}
k[u_1]\longrightarrow k[u_2],\quad p\longmapsto \theta p\theta^{-1}
\end{equation*}
gives an isomorphism.

Let $(G,V)$ be the $D_4$ case. Then $k[u_1]$ is a quadratic extension over $k$.
Hence there exist $a,b\in k$ with $b\not=0$ such that
$u_2=\theta (a+bu_1)\theta^{-1}$. Hence for
\begin{equation*}
g=\left(
	\theta, \theta^{-1}, \begin{pmatrix} 1&0\\ a&b\\\end{pmatrix}
\right)\in G_k,
\end{equation*}
we have $y=gx$.

Let $(G,V)$ be the $E_6$ case. 
There exists $p\in k[u_1]$ so that $u_2=\theta p\theta^{-1}$.
We claim that 
there exist $a,b,c,d\in k$ with  $ad-bc\not=0$ such that
$p=(c+du_1)/(a+bu_1)$.
In fact, if we consider the $k$-linear map %$\psi$ defined by
\begin{equation*}
\psi \colon k^4\longrightarrow k[u_1],\qquad (a,b,c,d)\longmapsto (a+bu_1)p-(c+du_1),
\end{equation*}
the kernel of $\psi$ is non-trivial.
Therefore there exists $(a,b,c,d)\in k^4\setminus \{0\}$
so that $(a+bu_1)p-(c+du_1)=0$, and since $p\notin k$,
we have $ad-bc\not=0$.

%Hence we have
%$$u_2=\theta\cdot  \frac{c+du_1}{a+bu_1}\cdot \theta^{-1}.$$
Hence for
\begin{equation*}
g=\left(
	\theta(a+bu_1)^{-1}, \theta^{-1}, \begin{pmatrix} a&b\\ c&d\\\end{pmatrix}
\right)\in G_k,
\end{equation*}
we have $y=gx$.
\end{proof}
We summarize the result in this subsection as follows.
\begin{prop}\label{rod-d4}
Let $(G,V)$ be the $D_4$ case.
\begin{enumerate}[{\rm (1)}]
\item
The map
\begin{equation*}
G_k\backslash V_k^\sst\longrightarrow \qsp(D),\quad
x\longmapsto \ti k(x)
\end{equation*}
is bijective.
\item
Let $x\in V_k^\sst$.
As a group over $k$, 
$G_x^\circ\cong (\gl(1)_{\ti k(x)})^2$.
%\item
%For any $x\in V_k^\sst$, $G_{xk}/G_{xk}^\circ\cong \Z/2\Z$.
\end{enumerate}
\end{prop}
\begin{prop}\label{rod-e6}
Let $(G,V)$ be the $E_6$ case.
\begin{enumerate}[{\rm (1)}]
\item
The map
\begin{equation*}
G_k\backslash V_k^\sst\longrightarrow \csp(D),\quad
x\longmapsto \ti k(x)
\end{equation*}
is bijective.
\item
Let $x\in V_k^\sst$.
As a group over $k$, 
$G_x^\circ\cong \gl(1)_{\ti k(x)}\times \gl(1)_k$.
%\item
%We have
%%
%\begin{equation*}
%G_{xk}/G_{xk}^\circ\cong
%\begin{cases}
%	\gS_3		&	k(x)=k		,\\
%	\Z/2\Z		&	[k(x):k]=2	,\\
%	\Z/3\Z		&	[k(x):k]=3	,\\
%	\{1\}		&	[k(x):k]=6	.
%\end{cases}
%\end{equation*}
%%
\end{enumerate}
\end{prop}

\subsection{Application to global fields}
If $k$ is a global field,
it is well known that
the sets $\qsp(D),\csp(D)$
can be described by means of local conditions.
Here, we review the argument.
We assume that $k$ is a global field in this subsection.
Also if $D$ is split, $\qsp(D)=\qsp$ and $\csp(D)=\csp$
for $n=2$ and $n=3$, respectively.
Hence we assume $D$ is non-split in this subsection.
Recall that $m=n^2$ is the dimension of $D$.

Let $\gM$ be the set of places of $k$.
For $v\in\gM$, let $k_v$ be the completion of $k$ at $v$.
We denote by $\inv_v(D)$ the Hasse invariant of $D\otimes k_v$ over $k_v$.
\begin{defn}
We define $\gM_D$
to be the set of elements $v\in\gM$ which satisfy $\inv_v(D)\not=0$.
\end{defn}
It is well known that $\gM_D$ is a finite set.
\begin{prop}
\begin{enumerate}[{\rm (1)}]
\item For $n=2$, the set $\qsp(D)$ consists of
elements $L\in\qsp$ such that
$L\otimes k_v$ is a quadratic extension of $k_v$
for all $v\in\gM_D$.
\item For $n=3$, the set $\csp(D)$ consists of
elements $L\in\csp$ such that
$L\otimes k_v$ is a cubic extension of $k_v$
for all $v\in\gM_D$.
\end{enumerate}
\end{prop}
\begin{proof}
We will prove the case $n=3$.
The case $n=2$ can be treated similarly.
%First note that if $D$ is split, the proposition is obvious
%since $\csp(D)=\csp$ and $\gM_D=\emptyset$.
%Hence we assume $D$ is non-split.
%Then any element of $\csp(D)$ is a cubic extension of $k$.

Let $L$ be an arbitrary separable cubic extension of $k$.
We denote by $\gM_L$ the set of places of $L$.
The field $L$ is an element of $\csp(D)$ if and only if
$D$ is split over $L$.
By the Hasse principle, this condition is equivalent to
that $D\otimes_k L_w\cong {\mathrm M(3,3)}_{L_w}$ for all $w\in\gM_L$.
Since $D\otimes_k k_v\cong {\mathrm M(3,3)}_{k_v}$ for all $v\notin\gM_D$,
we only need to consider $w$ which divides an element $v\in \gM_D$. 
For this $v$, $D_v=D\otimes k_v$ is a division algebra.
Hence for a separable extension $F/k_v$ with $[F: k_v]\leq 3$,
$D_v\otimes_{k_v}F\cong {\mathrm M(3,3)}_F$
if and only if $[F:k_v]=3$.
Therefore 
$D_v\otimes_{k_v}L_w\cong {\mathrm M(3,3)}_{L_w}$
if and only if $[L_w:k_v]=3$.
%Now the proposition holds
%because the restriction map between the Brauer groups
%of separable extension is the degree times map.
\end{proof}

\begin{rem}\label{exp}
Let $\ti T=\ker(G\rightarrow \gl(V))$.
Then it is easy to that

\begin{equation*}
\ti T=\{(t_{11},t_{12},t_2)\mid 
t_{11},t_{12},t_2\in\gl(1)_k, t_{11}t_{12}t_2=1\}\cong\gl(1)_k\times\gl(1)_k,
\end{equation*}
and hence
\begin{equation*}
G_x^\circ/\ti T\cong
\begin{cases}
(\gl(1)_{\ti k(x)}/\gl(1)_k)^2		& \text{the $D_4$ case},\\
\gl(1)_{\ti k(x)}/\gl(1)_k		& \text{the $E_6$ case}.\\
\end{cases}
\end{equation*}
For the non-split cases, $\ti k(x)$ is
a quadratic or cubic field over $k$
for any $x\in V_k^\sst$,
and hence $G_x^\circ/\ti T$
does not contain a split torus.
This shows that $(G,V)$ are of complete type for the non-split cases.

We conclude this subsection with a brief discussion
of the possible density theorems which we expect to derive
from the theory of the zeta function for our cases.

Roughly speaking, the global zeta function is a counting function
for the unnormalized Tamagawa numbers of $G_x^\circ/\ti T$ of points in
$x\in G_k\backslash V_k^\sst$.
Let $F=\ti k(x)$ and 
we denote by $h_F$ and $R_F$
the class number and  the regulator of $F$, respectively.
If we consider the canonical measure
on the adelization of $G_x^\circ/\ti T$,
the unnormalized Tamagawa number of this group is
$(\Res_{s=1}\zeta_F(s))^2$ (resp.\ $\Res_{s=1}\zeta_F(s)$)
for the $D_4$ case (resp.\ the $E_6$ case)
where $\zeta_F(s)$ is the Dedekind zeta function.
This leads us to believe that the theory of the zeta function
will eventually yield
the average density of $h_F^2R_F^2$ for $F\in\qsp(D)$ from the $D_4$ case, and
the average density of $h_FR_F$ for $F\in\csp(D)$ from the $E_6$ case.

So, we have to develop a local theory
to see that these are the correct interpretation of the problem.

\end{rem}

%\newpage
\section{The global zeta function}\label{gzf}
In this section, we study analytic properties
of the global zeta function for non-split cases.
The main result is Theorem \ref{ppf},
which describe the principal parts of the global zeta function.
\subsection{Preliminaries}\label{ss:pre}

%\ \linebreak[0]

In this subsection, we collect basic notations that we use in this section.
Also, we review some basic facts concerning adelic analysis
that we need later.
Throughout this section, $k$ is a number field.
Let $D$ be a non-split simple algebra of $k$ of dimension $4$ or $9$.
Then $D$ is a division algebra.
Since the argument is similar for the two cases,
we treat them simultaneously.
Recall that $m=n^2$ is the dimension of $D$.

Suppose that $G$ is a locally compact group and
$\Gamma$ a discrete subgroup of $G$
contained in the maximal unimodular subgroup of $G$.
For any left invariant measure $dg$ on $G$,
we choose a left invariant measure $dg$
(we use the same notation, but the meaning will be clear from the context)
on $X=G/\Gamma$ so that
$$\int_{G}f(g)\,dg=\int_X\sum_{\gamma\in\Gamma}f(g\gamma)\,dg.$$

Let $r_1$, $r_2$, $h_k$, $R_k$ and $\Delta_k$ be 
the number of real places, the number of complex
places, the class number, the regulator and 
the discriminant of $k$, respectively.
Let $e_k$ be the number of roots of unity contained in $k$.
We set
\begin{equation*}
\gc_k=2^{r_1}(2\pi)^{r_2}h_kR_ke_k^{-1}
.\end{equation*}

We refer \cite{weilc} as the basic reference for fundamental properties
on adeles.
%We assume that the reader is familiar with the basic definitions
%and facts concerning adeles. These may be found in \cite{weilc}.
The ring of adeles and the group of ideles
are denoted by $\A$ and $\ma$, respectively.
The adelic absolute value $|\;|$ on $\ma$ is normalized so that,
for $t\in\ma$, $|t|$ is the module of multiplication by $t$
with respect to any Haar measure $dx$ on $\A$, i.e. $|t|=d(tx)/dx$.
Let $\A^0=\{t\in\ma\mid |t|=1\}$.
We fix a non-trivial additive character $\ac{\ }$ of $\A/k$.
The set of positive real numbers is
denoted $\R_{+}$. 
Suppose $[k:\Q]=n$.
For $\lam\in\R_+$, $\lamb\in \ma$ 
is the idele whose component at any infinite place is $\lam^{1/n}$
and whose component at any finite place is $1$.  
Then we have $|\lamb|=\lam$.

We choose a Haar measure  $dx$ on $\A$ so 
that $\int_{\A/k} dx = 1$.
We define a Haar measure 
$\md t^0$ on $\A^0$ so that 
$\int_{\A^0/\mk} \md t^0 =  1$.  Using 
this measure, we choose a Haar
measure $\md t$ on $\ma$  so that 
$$
\int_{\ma} f(t)\, \md t  = \int_0^{\infty}   
\int_{\A^0}  f(\lamb t^0)\, \md \lam \md t^0
,$$ 
where $\md \lam = \lam^{-1}d\lam$ and $d\lam$ is a Lebesgue measure.

Let $\zeta_k(s)$ be the Dedekind zeta function of $k$.
We define
\begin{equation*}\label{Zkdefn}
Z_k(s)=|\Delta_k|^{s/2}
\left(\pi^{-s/2}\Gamma\left(\frac s 2\right)\right)^{r_1}
\left((2\pi)^{1-s}\Gamma(s)\right)^{r_2}\zeta_k(s)\,
.\end{equation*}
This definition differs from that in \cite{weilc}, p.129 by the
factor of $|\Delta_{k}|^{s/2}$ and from that in
\cite{yukiec} by the factor of $(2\pi)^{r_2}$. It is adopted here as
the most convenient for our purposes. It is well known that 
$\res_{s=1}Z_k(s)=\gc_k$. We define
\begin{equation*}
\phi(s)=\frac{Z_k(s)}{Z_k(s+1)},\quad
\varrho=\Res_{s=1}\phi(s)=\frac{\gc_k}{Z_k(2)},\quad
\text{and}\quad \gV_2=\varrho^{-1}.
\end{equation*}

For a complex variable $s$, we denote by $\Re(s)$ the real part.

For a vector space $V$ over $k$, $\av$ denotes its adelization.
Let $\sS(V_\A)$ be the spaces of 
Schwartz-Bruhat functions on $\av$.
We define the Haar measure $dx$ on $V_\A$
so that $\int_{V_\A/V_k}dx=1$.

We express elements of $G_2=\gl(2)$ as follows:
\begin{equation*} 
a(t_1,t_2)= \begin{pmatrix} t_1 & 0\\ 0 & t_2\end{pmatrix},\;
n(u) = \begin{pmatrix} 1 & 0\\ u & 1\end{pmatrix},\;
\nu=\begin{pmatrix} 0&1\\ 1&0\end{pmatrix}.
\end{equation*}

We recall the following well known facts concerning adelic analysis.
The proof may be found in \cite[Chapter 1]{yukiec}.
\begin{lem}\label{saa}
\begin{enumerate}[{\rm (1)}]
\item
Let $C\subset\gl(V)_\A$ be a compact set, and $\Phi\in\sS(V_\A)$.
Then there exists $\Psi\in\sS(\av)$ such that
$$|\Phi(gx)|\leq\Psi(x)$$ for all $g\in C, x\in V_\A$.
\item
Let $\Phi$ be a Schwartz-Bruhat function on $\A^n$.
Then there exist Schwartz-Bruhat functions $\Phi_1, \ldots, \Phi_n \geq 0$
such that $$|\Phi(x_1,\ldots,x_n)|\leq \Phi_1(x_1)\cdots\Phi_n(x_n).$$
\item
Suppose $\Phi\in\sS(\A)$. Then 
for any $N\geq1$,
$$\sum_{x\in k}\Phi(tx)\ll \max\{1,|t|^{-1}\}
,\qquad\sum_{x\in \mk}\Phi(tx)\ll |t|^{-N}.$$
\end{enumerate}
\end{lem}
\subsection{The global zeta function}\label{ss:gzf}

In this subsection, we define the global zeta function.
Also we give an estimate of an incomplete theta series
in order to use Shintani's lemma.

Recall that $G_2=\gl(2)$.
Let $T_2\subset G_2$ be the set of diagonal matrices and
$N_2\subset G_2$ be the set of lower-triangular matrices whose
diagonal entries are $1$. Then $B_2=T_2N_2$ is a Borel subgroup of
$G_2$. 

Let
\begin{align*}
G_{1\A}^0
&=	\{g_1=(g_{11},g_{12})\in G_{1\A}\mid |\nr(g_{11})|=|\nr(g_{12})|=1\},\\
G_{2\A}^0
&=	\{g_2\in G_{2\A}\mid |\det g_2|=1\},\\
\ag^0
&=	G_{1\A}^0\times G_{2\A}^0,\quad
		\overline G_\A=\R_+\times \ag^0,\\
\widehat T_{2\A}^0
&=	\{a(t_{21},t_{22})\mid t_{21},t_{22}\in\A^0\},\\	
T_{2\A}^0
&=	\{a(\mub^{-1},\mub)t_2\mid \mu\in\R_+, t_2\in \widehat T_{2\A}^0\},\\	
B_{2\A}^0
&=	T_{2\A}^0N_{2\A},\quad
		P_\A^0=G_{1\A}^0\times B_{2\A}^0.
\end{align*}
The group $\overline G_\A$ acts on $V_\A$ by assuming that
$\lam\in\R_+$ acts by multiplication by $\lamb$.
Throughout this section, we express elements
$\bar g\in \overline G_\A, g^0\in\ag^0$ as
\begin{equation*}
\bar g=(\lam,g_1,g_2), g^0=(g_1,g_2)
\end{equation*}
where $\lam\in\R_+,g_1\in G_{1\A}^0$, and $g_2\in G_{2\A}^0$.
We identify element $g^0\in \ag^0$ with $(1,g^0)\in\overline G_\A$
and $g_1\in G_{1\A}^0, g_2\in G_{2\A}^0$ with
$(1,g_1,1),(1,1,g_2)$.
We may also write as $\bar g=\lamb g^0$.

Let $\cK_2$ be the standard maximal compact subgroup of $G_{2\A}^0$ i.e.

\begin{equation*}
\cK_2
=	\prod_{v\in\gMr}{\mathrm O}(2,\R)\times
	\prod_{v\in\gMc}{\mathrm U}(2,\C)\times
	\prod_{v\in\gMf}\gl(2,\co_v)
\end{equation*}
Let $d\kappa _2$ be the Haar measure on $\cK_2$
such that the total volume of $\cK_2$ is $1$.

Let
\begin{equation*}
t_2=a_2(\mub^{-1}t_{21},\mub t_{22}),b_2=t_2n_2(u)
\end{equation*}
where $\mu\in \R_+,t_{21},t_{22}\in\A^0, u\in \A$.
Throughout this section, we assume that
\begin{equation*}
g_2=\kappa_2b_2=\kappa_2a_2(\mub^{-1}t_{21},\mub t_{22})n_2(u)
\end{equation*}
is the Iwasawa decomposition of $g_2\in G_{2\A}^0$.

The measure $du$ on $\mathbb{A}$ 
induces an invariant measure on $N_{\mathbb{A}}$. 
We put
\begin{equation*}
\md t=\md\mu\,\md t_{21}\,\md t_{22}, db_2=\mu^2\,\md t_2\,du_2
\end{equation*}
We use $dg_2=d\kappa_2db_2^0$ as the Haar measure on $G_{2\mathbb{A}}^0$.
It is well known that the volume of $G_{2\A}^0/G_{2k}$
with respect to the measure $dg_2$ is $\gV_2$.

We fix an arbitrary Haar measure $dg_1$ on $G_{1\A}^0$.
Since the rank of the group $G_1$ is $0$, $G_{1\A}^0/G_{1k}$ is compact.
We put
\begin{equation*}
\tau(G_1)=\int_{G_{1\A}^0/G_{1k}}dg_1.
\end{equation*}
We choose $dg^0=dg_1dg_2,d\bar g=\md\lam dg^0$
as Haar measures on $G_\A^0,\overline G_\A$, respectively.

For $\eta>0$, we define
\begin{equation*}
T_{2\eta+}^0=\{a(\mub^{-1},\mub)\mid \mu\in\R_+,\mu\leq \eta\}.
\end{equation*}
Let $\Omega_2\subset \widehat T_{2\A}^0N_{2\A}$
be a compact subset.
We define $\gS_2^0=\cK_2T_{2\eta+}^0\Omega_2$.
It is well known that for a suitable choice of $\eta$ and $\Omega_2$,
$\gS_2^0$ surjects to $G_{2\A}^0/G_{2k}$.
Also there exists another compact set $\widehat\Omega_2\in G_{2\A}^0$
such that $\gS_2^0\subset \widehat\Omega_2T_{2\eta+}^0$.
We fix a compact subset $\widehat\Omega_1\subset G_{1\A}^0$
which surjects to $G_{1\A}^0/G_{1k}$.
Let $\widehat\Omega=\widehat\Omega_1\times\widehat\Omega_2$.

\begin{defn}
Let $r\in\R$.
We define $C(\ag^0/\rg,r)$ to be
the set of continuous functions $f(g^0)$ on $\ag^0/\rg$
satisfying
\begin{equation*}
\sup_{g^0\in \widehat\Omega T_{2\eta+}^0}f(g^0)\mu^{-r}<\infty.
\end{equation*}
A function $f$ on $\ag^0/\rg$ is said to be {\em slowly increasing}
if $f\in C(\ag^0/\rg,r)$ for some $r\in\R$.
\end{defn}

Note that $C(\ag^0/\rg,r_1)\subset C(\ag^0/\rg,r_2)$ if $r_1>r_2$
and $C(\ag^0/\rg,r)\subset L^1(\ag^0/\rg,dg^0)$ if $r>-2$.

\begin{defn}
Let $\Phi\in\sS(V_\A)$ and $\bar g\in \overline G_\A$.
For any subset $L\subset V_k$, we define
an incomplete theta series $\Theta_L(\Phi,\bar g)$ by
\begin{equation*}
\Theta_L(\Phi,\bar g)=\sum_{x\in L} \Phi(\bar gx).
\end{equation*}
If $L$ is invariant under the action of a subgroup $H_k\subseteq G_k$,
we often regard $\Theta_L(\Phi,\bar g)$ as a function on $\bar \ag/H_k$.
For example, if $L=V_k^\sst$, $\Theta_{V_k^\sst}(\Phi,\bar g)$
is a function on $\overline G_\A/G_k$.
\end{defn}

\begin{lem}\label{Theta}
For any $N\geq1$,
\begin{equation*}
\Theta_{V_k^\sst}(\Phi,\bar g)\ll
\begin{cases}
\lam^{-2N-2}\mu^{-1}	&	\lam\geq1,\\
\lam^{-2-m}\mu^{-1}		&	\lam\leq1,
\end{cases}
\end{equation*}
for $\bar g\in\R_+\times\gS^0$. 
\end{lem}
\begin{proof}
By (1) of Lemma \ref{saa}, we may assume
$\bar g=\lamb(1,a_2(\mub^{-1},\mub)),\mu\ll 1$.
For $x=(x_1,x_2)\in V_k^\sst$, we have $x_1\not=0$ and $x_2\not=0$.
Note that the weight of $a(t_1,t_2)\in T_2$ with respect to
each $k$-coordinate of $x_1$ and $x_2$ is $t_1$ and $t_2$, respectively.
Hence, by (2) and (3) of Lemma \ref{saa},
for any $N_1,N_2\geq1$,
\begin{align*}
\Theta_{V_k^\sst}(\Phi,\bar g)
&\ll	(\lam^{-1}\mu)^{N_1}(\lam^{-1}\mu^{-1})^{N_2}
			\max(1,\lam^{-1}\mu)^{m-1}
			\max(1,\lam^{-1}\mu^{-1})^{m-1}\\
&\leq	\lam^{-N_1-N_2}\mu^{N_1-N_2}
			\max(1,\lam^{2-2m})\max(1,\mu^{m-1})\max(1,\mu^{1-m})\\
&\ll		\lam^{-N_1-N_2}\max(1,\lam^{2-2m})\mu^{N_1-N_2+1-m}.
\end{align*}
Note that $\max\{1,ab\}\leq\max\{1,a\}\cdot\max\{1,b\}$ for $a,b\geq0$.
For $\lam\geq1$, take $N_1=N+m-2,N_2=N$.
For $\lam\leq1$, take $N_1=3,N_2=1$.
Then we have the proposition.
\end{proof}

Now we define the global zeta function.
\begin{defn}
For $\Psi\in\sS(V_\A)$ and a complex variable $s$,
we define
\begin{align*}
Z(\Phi,s)
&=	\int_{\overline G_\A/\rg}
		\lam^s\Theta_{V_k^\sst}(\Phi,\bar g)\, d\bar g,\\
Z_+(\Phi,s)
&=	\int_{\underset{\lam\geq1}{\overline G_\A/\rg}}
		\lam^s\Theta_{V_k^\sst}(\Phi,\bar g)\, d\bar g.
\end{align*}
\end{defn}
The integral $Z(\Phi,s)$ is called the {\em global zeta function}.
%\fbox{comment} on H. Saito's work.
By Lemma \ref{Theta}, the integral $Z(\Phi,s)$ converges absolutely
and locally uniformly on a certain right half-plane and the integral
$Z_+(\Phi,s)$ is an entire function.
Since the global zeta function is well defined for $V_k^\sst$,
%This means that $(G,V)$ is of complete type and therefore,
by Theorem (0.3.7) in \cite{yukiec} (which is due to Shintani),
$Z(\Phi,s)$ can be continued meromorphically to the entire plane and
satisfies a functional equation
\begin{equation*}
Z(\Phi,s)=Z(\widehat\Phi,2m-s).
\end{equation*}
The purpose of this section is to determine the pole structure
and to describe the residues by means of certain distributions.
\begin{rem}
The above definition of the zeta function looks slightly different
from the original definition in \cite{sash}, but
these are essentially the same.
%For the convenience of the reader,
We briefly compare these functions.
Let $\ti G=G/\ti T$. Recall that we put $\ti T=\ker(G\rightarrow \gl(V))$.
Let $d\ti g$ be an invariant measure on $\ti G_\A$.
The original definition of the global zeta function is as follows:
\begin{equation*}
Z^\ast(\Phi,s)
=	\int_{\ti G_\A/\ti G_k}
		|\chi(\ti g)|^s\Theta_{V_k^\sst}(\Phi,\ti g)\, d\ti g.
\end{equation*}
Since $\ti T\cong\gl(1)\times\gl(1)$ is a split torus,
the first Galois cohomology set $H^1(k',\ti T)$ is trivial
for any field $k'$ containing $k$.
This implies that the set of $k'$-rational point of $\ti G$
coincides with $G_{k'}/\ti T_{k'}$.
Therefore $\ti G_\A=G_\A/\ti T_\A$ and
$\ti G_\A/\ti G_k=G_\A/\ti T_\A G_k$.
Let $\ti T_\A^0=\ag^0\cap \ti T_\A$.
Then we have
\begin{equation*}
(\R_+\times G_\A^0)/\ti T_\A^0\cong G_\A/\ti T_\A
\end{equation*}
via the map which sends the class of
$(\lam,g_{11},g_{12},g_2)$ to class of $(g_{11},g_{12},\lamb g_2)$.
Moreover, this map is
compatible with their actions on $V_\A$.
If we identify $\overline G_\A/\ti T_\A^0$ with $\ti G_\A$ via the isomorphism,
then we have
$|\chi(\bar g)|=\lam^4$ for the $D_4$ case and
$|\chi(\bar g)|=\lam^{12}$ for the $E_6$ case.
Also the volume of $\ti T_\A^0/\ti T_k\cong (\A^0/\mk)^2$ is finite.
Hence $Z(\Phi,4s)$ is a constant multiple of $Z^\ast(\Phi,s)$
for the $D_4$ case and $Z(\Phi,12s)$ is a constant multiple
of $Z^\ast(\Phi,s)$ for the $E_6$ case.
(The constant depends on the choice of the measure.)
Our choice of $Z(\Phi,s)$ is 
for the conventions of our global theory.
%To deduce density theorems,
%we will probably use $Z(\Phi,s)$.
\end{rem}

Let
\begin{equation*}
{\mathrm M}\Phi(x)=\int_{\cK_2} \Phi(\kappa_2 x)\, d\kappa_2.
\end{equation*}
Then $Z(\Phi,s)=Z({\mathrm M}\Phi,s)$ and
${\mathrm M}\Phi(x)$ is $\cK_2$-invariant.
Therefore, we may assume the following for the rest of this section.
\begin{asmp}
The Schwartz-Bruhat function $\Phi$ is $\cK_2$-invariant.
\end{asmp}

\subsection{The principal part}\label{ss:pp}

For $x=(x_1,x_2)$ and $y=(y_1,y_2)$, we define
\begin{equation*}
[x,y]=\trc(x_1y_2)+\trc(x_2y_1).
\end{equation*}
This is a non-degenerate bilinear form on $V$.
For $\bar g=(\lam,g_{11},g_{12},g_2)\in\R_+\times G_\A^0$,
we define
\begin{equation*}
\bar g^\iota=(\lam^{-1},g_{12}^{-1},g_{11}^{-1},\nu{}^tg_2^{-1}\nu).
\end{equation*}
This is an involution and the above bilinear form satisfies
\begin{equation*}
[\bar g x,\bar g^\iota y]=[x,y].
\end{equation*}
Recall that $\ac{\ }$ is a non-trivial additive character of $\A/k$.
For $\Phi\in\sS(V_\A)$, we define its Fourier transform by
\begin{equation*}
\widehat\Phi(x)=\int_{V_\A}\Phi(y)\ac{[x,y]}\, dy.
\end{equation*}
It is easy to see that the Fourier transform of $\Phi(\bar g\cdot)$ is
$\lam^{-8}\widehat\Phi(\bar g^\iota\cdot)$.

For $\lam\in\R_+$, we define $\Phi_\lam(x)=\Phi(\lamb x)$.
\begin{defn}
For $\Phi\in\sS(V_\A),s\in \C$ and $g^0\in G_\A^0$, we define
\begin{align*}
J(\Phi,g^0)
&=	\sum_{x\in S_k}\widehat\Phi((g^0)^\iota x)-\sum_{x\in S_k}\Phi(g^0x),\\
I^0(\Phi)
&	=\int_{G_\A^0/G_k}J(\Phi,g^0)\, dg^0,\\
I(\Phi,s)
&	=\int_0^1\lam^s I^0(\Phi_\lam)\, \md\lam.
\end{align*}
\end{defn}
Then by the Poisson summation formula, we have the following.
\begin{prop}
We have
\begin{equation*}
Z(\Phi,s)=Z_+(\Phi,s)+Z_+(\widehat\Phi,2m-s)+I(\Phi,s).
\end{equation*}
\end{prop}
We study the last term for the rest of this section.

\subsection{Stratification}\label{ss:str}
In this subsection, we consider a stratification of $V_k$.
Let
\begin{equation*}
Y_1=\{x\in V\mid x_1=0\},\quad Y_1^\sst=\{x\in Y_1\mid x_2\not=0\}.
\end{equation*}
We define $S_1=GY_1^\sst$. Let $P=G_1\times B_2$.
\begin{lem}\label{stratification}
We have
\begin{enumerate}[{\rm (1)}]
\item
$V_k\setminus\{0\}=V_k^\sst\amalg S_{1k}$,
\item
$S_{1k}=G_k\times_{P_k}Y_{1k}^\sst$.
\end{enumerate}
\end{lem}
\begin{proof}
We consider (1).
Let $x\in V_k\setminus\{0\}$ and $x\notin V_k^\sst$.
Since either $x_1\not=0$ or $x_2\not=0$,
there exists an element $g\in G_k$
such that the first coordinate of $gx$ is $1$.
Replacing $x$ by $gx$, we may assume that
$x$ is of the form $x=(1,-x_2)$, where $x_2\in W_k$.
Then $F_x(v)$ is the characteristic polynomial of $x_2$
and the condition $P(x)=0$ is equivalent to that
the characteristic polynomial of $x_2$ has a multiple root.

Let $L=k[x_2]$ be the subalgebra of $D_k$ generated by $x_2$ over $k$.
Since $D_k$ has no zero divisor, $k[x_2]$ is a (commutative) integral domain
which is finite over the field $k$. So it is a field.
Then since the degree of extension $[L:k]$ divides $\dim_kD_k=n^2$, 
it is either $1$ or $n$.
Note that we are assuming $n=2$ or $3$.
Assume $[L:k]=n$. Then $F_x(v)$ is a minimum polynomial of $x_2$ over $k$
because the degree of $F_x(v)$ is $n$.
Since any field extension of an algebraic number field is separable,
we conclude that $F_x(v)$ does not have a multiple root.
This is an contradiction and hence $[L:k]=1$,
which implies $x_2\in k$.
Therefore, there exists an element $g_2\in G_{2k}$
such that $g_2x\in Y_{1k}^\sst$.
This proves (1).

It is easy to see that $P_{1k}Y_{1k}^\sst=Y_{1k}^\sst$
and that if $x\in Y_{1k}^\sst, g\in G_k$ and $gx\in Y_{1k}^\sst$
then $g\in P_k$. This proves (2).
\end{proof}

\subsection{The smoothed Eisenstein series}\label{ss:ses}
To compute $I^0(\Phi)$,
it seems natural to divide the
index set $S_k$ of the summation
into its $G_k$-orbits
and perform integration separately.
However, we can not put this into practice
because the corresponding integrals diverge.
This is the main difficulty
when one calculates the global zeta functions
of the prehomogeneous vector spaces.
To surmount this problem
Shintani \cite{shintania} introduced
the smoothed Eisenstein series of $\gl(2)$.
He used this series to determine the principal parts
of the global zeta functions for the space of binary cubic forms
and the space of binary quadratic forms.
Later A. Yukie \cite{yukiec} generalized
the theory of Eisenstein series
to the products of $\gl(n)$'s,
and applied it to determine the principal parts of
the global zeta functions in some cases.
In this subsection, we essentially repeat
the argument of Shintani and Yukie in our settings.
%In this subsection, we define distributions
%related to the smoothed Eisenstein series.
%This is a way to treat the divergent integrals
%and first introduced by Shintani \cite{shintania}.

We express the Iwasawa decomposition of $g_2\in G_{2\A}^0$ as
\begin{equation*}
g_2=\kappa_2(g_2)a_2(t_{21}(g_2),t_{22}(g_2))n_2(u(g_2)).
\end{equation*}
Let $s\in\C$.
The Eisenstein series of $G_{2\A}^0$ for $B_2$ is defined as
\begin{equation*}
E(g_2,s)=\sum_{\gamma\in G_{2k}/B_{2k}}|t_{21}(g_2\gamma)|^{s+1}
\end{equation*}
It is well known that the summation defining $E(g_2,s)$
converges absolutely and locally uniformly in certain right half-plane
and can be continued meromorphically to the whole complex plane.
For analytic properties of $E(g_2,s)$, see \cite{wright}, \cite{yukiec}.

Let $\psi(s)$ be an entire function of $s$ such that
\begin{equation*}
\sup_{c_1<\Re(s)<c_2}(1+|s|^N)|\psi(s)|<\infty
\end{equation*}
for all $c_1<c_2, N>0$.
Moreover, we assume $\psi(1)\not=0$.
For a complex variable $w$, we define
\begin{equation*}
\Lambda_\psi(w;s)=\frac{\psi(s)}{w-s}.
\end{equation*}
Where there is no confusion, we drop $\psi$ and use the notation
$\Lam(w;s)$ instead.
\begin{defn}
For a complex variable $w$, we define
\begin{equation*}
\sE(g^0,w,\psi)=\ct\int_{\Re(s)=r_1}E(g_2,s)\Lam_\psi(w;s)\, ds.
\end{equation*}
for some $r_1>1$.
\end{defn}
Note that the above definition does not depend on the choice of $r_1$.
The function $\sE(g^0,w,\psi)$ is called the smoothed Eisenstein series.
When there is no confusion, we drop $\psi$ and use the notation
$\sE(g^0,w)$ instead.

The following proposition is known as Shintani's lemma.

\begin{prop}\label{ses}
\begin{enumerate}[{\rm (1)}]
\item
The function $\sE(g^0,w)$ is holomorphic for $\Re(w)>0$
except for a simple pole at $w=1$ with the residue $\varrho\psi(1)$.
\item
Let $f\in C(G_\A^0/G_k,r)$ for some $r>-2$.
Then the integral
\begin{equation*}
\int_{\ag^0/\rg}f(g^0)\sE(g^0,w)\, dg^0
\end{equation*}
becomes a holomorphic function for $\Re(w)>1-\epsilon$
for a constant $\epsilon>0$
except possibly for a simple pole at $w=1$ with residue
\begin{equation*}
\varrho\psi(1)\int_{\ag^0/\rg}f(g^0)\, dg^0.
\end{equation*}
\item
For a slowly increasing function $f(g^0)$ on $\ag^0/\rg$,
the integral
\begin{equation*}
\int_{\ag^0/\rg}f(g^0)\sE(g^0,w)\, dg^0
\end{equation*}
becomes a holomorphic function on a certain right half-plane.
\item
We have
\begin{equation*}
\int_{\ag^0/\rg}\sE(g^0,w)\, dg^0=\tau(G_1)\Lam(w;1).
\end{equation*}
\end{enumerate}
\end{prop}
The above proposition was first proved for $\gl(2)$ by Shintani
\cite[pp. 172, 173, 177]{shintania}.
The adelic proof is given in \cite[pp. 527, 528]{wright}.
In our case, we included the first factor $G_{1\A}^0$
in the statement instead of just considering $\gl(2)$,
but exactly the same proof works because $G_{1\A}^0/G_{1k}$ is compact.
For the convenience of the reader, we give a sketch of the proof here.
To prove Proposition \ref{ses}, the following lemma 
on the Eisenstein series for $G_2=\gl(2)$ is crucial.
For the proof, see \cite{shintania}, \cite{wright}, or \cite{yukiec}.
\begin{lem}\label{es}
\begin{enumerate}[{\rm (1)}]
\item
Let $E_{N_2}(g_2,s)$ be the constant term of
$E(g_2,w)$ with respect to $N_2$ i.e.
$$E_{N_2}(g_2,s)=\int_{N_{2\A}/N_{2k}}E(g_2n_2(u),s)\ du.$$
The constant term has the following explicit formula:
\begin{equation*}
E_{N_2}(g_2,s)=\mu^{-s-1}+\mu^{s-1}\phi(s).
\end{equation*}
\item
Let $\ti E(g_2,s)=E(g_2,s)-E_{N_2}(g_2,s)$ be the non-constant term.
Then $\ti E(g_2,s)$ is holomorphic for $\Re(s)>0$.
Moreover, for any $s$ in this region and $l>1$,
\begin{equation*}
|\ti E(g_2,s)|\ll \mu^{2l-1}.
\end{equation*}
\item
We have
\begin{equation*}
\int_{G_{2\A}^0/G_{2k}}\sE(g_2,w)\, dg_2=\Lambda(w,1).
\end{equation*}
\end{enumerate}
\end{lem}
%Roughly speaking, this proposition says that
%we can multiply $\sE(g^0,w)$ to any slowly increasing
%function and make it integrable if $\Re(w)\gg0$.
%Moreover, if the function is slightly better than integrable,
%the resulting integral, as a function of $w$,
%has a simple pole at $w=1$ with the residue a constant multiple
%of the original integral.
%This is the key idea to separate contributions from strata.
\noindent
{\it Sketch of the proof of Proposition \ref{ses}}.
Let 
\begin{equation*}
E'(g_2,s)=E(g_2,s)-\mu^{s-1}\phi(s)=\mu^{-s-1}+\ti E(g_2,s).
\end{equation*}
This function is holomorphic for $\Re(s)>0$.
Recall that $\rho=\Res_{s=1}\phi(s)$.

Let $\epsilon$ be a small positive number.
By moving the contour, we have
\begin{align*}
\ct\int_{\Re(s)=r_1}\mu^{s-1}\phi(s)\frac{\psi(s)}{w-s}\, ds
&=	\frac{\rho\psi(1)}{w-1}
	+\ct\int_{\Re(s)=\epsilon}\mu^{s-1}\phi(s)\frac{\psi(s)}{w-s}\, ds,\\
\ct\int_{\Re(s)=r_1}E'(g_2,s)\frac{\psi(s)}{w-s}\, ds
&=	\ct\int_{\Re(s)=\epsilon}E'(g_2,s)\frac{\psi(s)}{w-s}\, ds,
\end{align*}
and each of the integrals in the right hand sides
is holomorphic for $\Re(w)>\epsilon$.
This proves (1).

We consider (2).
Since $C(G_\A^0/G_k,r_1)\subset C(G_\A^0/G_k,r_2)$ for $r_1>r_2$,
we may assume $-2<r<0$. By moving the contour, we have
\begin{align*}
&\ct\int_{\Re(s)=r_1}
	\left(\int_{\ag^0/\rg}f(g^0)\mu^{s-1}\, dg^0\right)
		\phi(s)\frac{\psi(s)}{w-s}\, ds\\
&\quad=\frac{\rho\psi(1)}{w-1}\int_{\ag^0/\rg}f(g^0)\, dg^0\\
&\qquad	+\ct\int_{\Re(s)=-r/2}
		\left(\int_{\ag^0/\rg}f(g^0)\mu^{s-1}\, dg^0\right)
			\phi(s)\frac{\psi(s)}{w-s}\, ds.
\end{align*}
The second term in the right hand side is holomorphic if $\Re(w)>-r/2$
because for $\Re(s)=-r/2$, we have
$f(g^0)\mu^{s-1}\in C(\ag^0/\rg,(r-2)/2)\subset L^1(\ag^0/\rg,dg^0)$.
A similar argument shows that
\begin{equation*}
\ct\int_{\Re(s)=r_1}
	\left(\int_{\ag^0/\rg}f(g^0)E'(g_2,s)\, dg^0\right)
		\frac{\psi(s)}{w-s}\, ds
\end{equation*}
is holomorphic for $\Re(w)>0$. This proves (2).

To see (3), we may assume
$f(g^0)\in C(\ag^0/\rg,r)$ for some $r\leq-2$.
Then $f(g^0)\mu^{s-1}\in C(\ag^0/\rg,-1)$ for $\Re(s)\geq -r$,
$f(g^0)\mu^{-s-1}\in C(\ag^0/\rg,-1)$ for $\Re(s)\leq r$,
and $f(g^0)\ti E(g_2,s)\in C(\ag^0/\rg,-1)$ for $\Re(s)>\epsilon$
where $\epsilon>0$ is a small number.
Hence by moving the contour for each of these terms, we obtain (3).
Now (4) immediately follows from Lemma \ref{es} (3),
and these finish the proof of Proposition \ref{ses}.
\hspace*{\fill}$\square$

%\vspace{1zw}

Let $\sE_{N_2}(g^0,w)$ be the constant term of $\sE(g^0,w)$
with respect to $N_2$ i.e.
\begin{equation*}
\sE_{N_2}(g^0,w)=\int_{N_{2\A}/N_{2k}}\sE(g^0n(u),w)\, du.
\end{equation*}
By Lemma \ref{es} (1), we have
\begin{equation*}
\sE_{N_2}(g^0,w)=\ct\int_{\Re(s)=r_1}(\mu^{-s-1}+\mu^{s-1}\phi(s))\Lam(w;s)\, ds.
\end{equation*}
\begin{defn}
Let $f(w),g(w)$ be holomorphic functions of $w\in\C$ in some right half-plane.
We use the notation $f(w)\sim g(w)$ if $f(w)-g(w)$ can be continued
meromorphically to $\{w\mid \Re(w)>1-\epsilon\}$ for some $\epsilon >0$
and is holomorphic at $w=1$.
\end{defn}

We define
\begin{equation*}
I^0(\Phi,w)=\int_{\ag^0/\rg}J(\Phi,g^0)\sE(g^0,w)\, dg^0.
\end{equation*}
By Lemma \ref{Theta}, $J(\Phi,g)\in C(\ag^0/\rg,-1)$.
Hence, by Proposition \ref{ses} (2), we have the following.
\begin{prop}\label{I^0}
We have
\begin{equation*}
I^0(\Phi,w)\sim \varrho\Lam(w;1)I^0(\Phi).
\end{equation*}
\end{prop}
\begin{defn}
For a complex variable $w$, we define
\begin{align*}
\Xi_1(\Phi,w)&=\int_{\ag^0/\rg}\Theta_{S_{1k}}(\Phi,g^0)\sE(g^0,w)\, dg^0,\\
\Xi_\#(\Phi,w)&=\Phi(0)\int_{\ag^0/\rg}\sE(g^0,w)\, dg^0.
\end{align*}
\end{defn}
Since $\Theta_{S_k'}(\Phi,g^0)$ is a slowly increasing function,
by Proposition \ref{ses} (3), the integral $\Xi_1(\Phi,w)$ converges absolutely
for sufficiently large $\Re(w)$.
It is proved in \cite{yukiec} that $E(g_2,s)=E({}^tg_2^{-1},s)$
for $g_2\in G_{2\A}^0$.
Hence, $\sE(g^0,w)=\sE((g^0)^\iota,w)$ for $g^0\in \ag^0$.
Therefore, by Lemma \ref{stratification}, we have the following.
\begin{prop}\label{Idecomposition}
We have
\begin{equation*}
I^0(\Phi,w)
=	\Xi_1(\widehat\Phi,w)+\Xi_\#(\widehat\Phi,w)
	-\Xi_1(\Phi,w)-\Xi_\#(\Phi,w).
\end{equation*}
\end{prop}
For $\Xi_\#(\Phi,w)$, Proposition \ref{ses} immediately leads
to the following.
\begin{prop}\label{xisharp}
We have
\begin{equation*}
\Xi_\#(\Phi,w)=\Phi(0)\tau(G_1)\gV_2\varrho\Lam(w,\rho).
\end{equation*}
\end{prop}
We study $\Xi_1(\Phi,w)$ in \S \ref{ss:cus}.

\subsection{The zeta function associated with the space of division algebra}
\label{ss:zsa}
Since the prehomogeneous vector space $(G_1,W)$ of (single) division algebra
appears in the induction process, we have to know the principal part of the
zeta function for this case.
This function is essentially the same as that of Godement-Jacquet \cite{goja}.
In this subsection we review the principal part
of the zeta function in this case.

We put $P_1(x)=\nr(x)$ for $x\in W$ and 
$W^\sst=\{x\in W\mid P_1(x)\not=0\}$.
Note that $W_k^\sst=\{x\in W_k\mid x\not=0\}$.
The group $\R_+\times G_{1\A}^0$ acts on $W_\A$ by assuming that
$\lam\in\R_+$ acts by multiplication by $\lamb$.
For any subset $L\subset W_k$,
we define $\Theta_L(\Psi,\lamb g_1)$ in the obvious manner.

\begin{defn}
For $\Psi\in\sS(W_\A)$ and $s\in\C$,
\begin{align*}
Z_W(\Psi,s)
&=	\int_{\R_+\times G_{1\A}^0/G_{1k}}
		\lam^s\Theta_{W_k^\sst}(\Psi,\lamb g_1)\, \md\lam dg_1^0,\\
Z_{W+}(\Psi,s)
&=	\int_{\underset{\lam\geq1}{\R_+}\times G_{1\A}^0/G_{1k}}
		\lam^s\Theta_{W_k^\sst}(\Psi,\lamb g_1)\, \md\lam dg_1^0.
\end{align*}
\end{defn}
The following lemma is a direct consequence of Lemma \ref{saa}.
\begin{lem}
The integral defining $Z_W(\Phi,s)$ converges absolutely
and locally uniformly in the region $\Re(s)>m$, and
the integral defining $Z_{W+}(\Phi,s)$ 
is an entire function.
\end{lem}
For $x,y\in W$, we put
\begin{equation*}
[x,y]_W=\trc(xy)
\end{equation*}
This defines a non-degenerate bilinear form on $W$.
We note that this bilinear form satisfies
$[g_1x,g_1^\iota y]_W=[x,y]_W$ where we put
$(g_{11},g_{12})^\iota=(g_{12}^{-1},g_{11}^{-1})$.

We define the Fourier transform on $\sS(W_\A)$ by
\begin{equation*}
\Psi^\ast(x)=\int_{W_\A}\Psi(y)\ac{[x,y]_W}\, dy.
\end{equation*}
Then by the \psf, we have
\begin{equation*}
\Theta_{W_k^\sst}(\Psi,\lamb g_1)
=	\lam^{-m}\Theta_{W_k^\sst}(\Psi^\ast,\lamb^{-1} (g_1)^\iota)
	+\lam^{-m}\Psi^\ast(0)-\Psi(0).
\end{equation*}
Applying the above equation,
 we obtain the following principal part formula for this zeta function.
\begin{prop}\label{Z_W}
We have
\begin{equation*}
Z_W(\Psi,s)
=	Z_{W+}(\Psi,s)+Z_{W+}(\Psi^\ast,m-s)
	+\tau(G_1)\left(\frac{\Psi^\ast(0)}{s-m}-\frac{\Psi(0)}{s}\right),
\end{equation*}
where $Z_{W+}(\Psi,s)$ and $Z_{W+}(\Psi^\ast,m-s)$ are entire functions.
\end{prop}

\subsection{Contribution from unstable strata}\label{ss:cus}
In this subsection, we express the residue of $\Xi_1(\Phi,w)$
in terms of $Z_W$ defined in the previous subsection.
We identify $Y_1$ (see \S \ref{ss:str}) with the space $W$
of single division algebras in \S \ref{ss:zsa}.

\begin{defn}
For $\Phi\in\sS(V_\A)$,
we define a Schwartz-Bruhat function $\cR_W\Phi$ on
$W_\A$ by restricting to $Y_{1\A}$.
\end{defn}
\begin{prop}\label{xi1}
By changing $\psi$ if necessary, we have
\begin{equation*}
\Xi_1(\Phi,w)\sim \rho\Lam(w;1)Z_W(\cR_W\Phi,2).
\end{equation*}
\end{prop}
\begin{proof}
We have
\begin{align*}
\Xi_1(\Phi,w)
&=	\int_{\ag^0/\rg}\Theta_{S_{1k}}(\Phi,g^0)\sE(g^0,w)\, dg^0\\
&=	\int_{\ag^0/P_k}\Theta_{Y_{1k}^\sst}(\Phi,g^0)\sE(g^0,w)\, dg^0\\
&=	\int_{P_\A^0/P_k}\Theta_{Y_{1k}^\sst}(\Phi,p^0)\sE(p^0,w)\, dg^0\\
&=	\int_{G_{1\A}^0/G_{1k}\times B_{2\A}^0/B_{2k}}
\Theta_{Y_{1k}^\sst}(\Phi,(g_1,b_2))\sE(b_2,w)\, dg_1db_2\\
&=	\int_{G_{1\A}^0/G_{1k}\times T_{2\A}^0/T_{2k}}
\mu^2\Theta_{Y_{1k}^\sst}(\Phi,(g_1,t_2))\sE_{N_2}(t_2,w)\, dg_1\md t_2.
\end{align*}
The last step is because $N_2$ acts on $Y_1$ trivially.
By changing $g_{11}$ to $g_{11}t_{22}^{-1}$, we have
\begin{equation*}
\Xi_1(\Phi,w)
=	\int_{\R_+\times G_{1\A}^0/G_{1k}}\mu^2\Theta_{W_k^\sst}(\cR_W\Phi,(\mu,g_1))
	\sE_{N_2}(a(\mub^{-1},\mub),w)\, \md\mu dg_1.
\end{equation*}
By the definition of $Z_W(\Psi,s)$, we have
\begin{equation*}
\int_{\R_+\times G_{1\A}^0/G_{1k}}
	\mu^{\mp s+1}\Theta_{W_k^\sst}(\cR_W\Phi,(\mu,g_1))\, \md\mu dg_1
=Z_W(\cR_W\Phi,\mp s+1)
\end{equation*}
for $\Re(s)<-3,\Re(s)>3$, respectively.
Since
\begin{equation*}
\ct\int_{\Re(s)=r_2<-3}Z_W(\cR_W\Phi,-s+1)\Lam(w;s)\, ds\sim 0,
\end{equation*}
we have
\begin{equation*}
\Xi_1(\Phi,w)
\sim
\ct\int_{\Re(s)=r_3>3}Z_W(\cR_W\Phi,s+1)\phi(s)\Lam(w;s)\, ds.
\end{equation*}
By Proposition \ref{Z_W}, $Z_W(\cR_W\Phi,s+1)$ has a simple pole
at $s=3$ and is holomorphic for $\Re(s)>0$ except for that point.
If we consider $(s-3)\psi(s)$ instead of $\psi(s)$,
this function still satisfies the property we have assumed. 
Namely,
\begin{equation*}
\sup_{c_1<\Re(s)<c_2}(1+|s|^N)|(s-3)\psi(s)|<\infty
\end{equation*}
for all $c_1<c_2, N>0$ and $(s-3)\psi(s)|_{s=1}\not=0$.
Therefore, by changing $\psi(s)$ to $(s-3)\psi(s)$,
we may assume that
$Z_W(\cR_W\Phi,s+1)\Lam(w;s)$
is holomorphic for $\Re(s)>0$.
(This is the passing principle (3.6.1) of \cite{yukiec}.)
Hence,
\begin{align*}
\Xi_1(\Phi,w)
&\sim
\ct\int_{\Re(s)=1/2}Z_W(\cR_W\Phi,s+1)\phi(s)\Lam(w;s)\, ds\\
&\quad+\varrho\Lam(w;1)Z_W(\cR_W\Phi,2)\\
&\sim\varrho\Lam(w;1)Z_W(\cR_W\Phi,2)
\end{align*}
This proves the proposition.
\end{proof}
\subsection{The principal part formula}\label{ss:ppf}

\begin{thm}\label{ppf}
Suppose that $\Phi={\mathrm M}\Phi$. Then
\begin{align*}
Z(\Phi,s)
&=Z_+(\Phi,s)+Z_+(\widehat\Phi,2m-s)\\
&\quad
	+\tau(G_1)\gV_2\left(\frac{\widehat\Phi(0)}{s-2m}-\frac{\Phi(0)}{s}\right)
	+\frac{Z_W(\cR_W\widehat\Phi,2)}{s-(2m-2)}-\frac{Z_W(\cR_W\Phi,2)}{s-2},
\end{align*}
where the first two terms in the right hand side are entire functions.
\end{thm}
\begin{proof}
By Proposition \ref{Idecomposition}, \ref{xisharp} and \ref{xi1},
\begin{equation*}
I^0(\Phi,w)=\rho\Lam(w;1)
	\left(
		\tau(G_1)\gV_2(\widehat\Phi(0)-\Phi(0))
		+Z_W(\cR_W\Phi,2)-Z_W(\cR_W\widehat\Phi,2)
	\right)
\end{equation*}
for a suitable choice of $\psi(s)$.
Hence, together with Proposition \ref{I^0}, we obtain
\begin{equation*}
I^0(\Phi)
=		\tau(G_1)\gV_2(\widehat\Phi(0)-\Phi(0))
		+Z_W(\cR_W\Phi,2)-Z_W(\cR_W\widehat\Phi,2).
\end{equation*}
Recall that $I(\Phi,s)=\int_0^1 I^0(\Phi_\lam)\md\lam$ where
$\Phi_\lam(x)=\Phi(\lamb x)$.
It is easy to see that
\begin{equation*}
\Phi_\lam(0)=\Phi(0),\quad \widehat{\Phi_\lam}(0)=\lam^{-2m}\widehat\Phi(0).
\end{equation*}
Since
\begin{align*}
Z_W(\cR_W\Phi_\lam,s)&=\lam^{-s}Z_W(\cR_W\Phi,s),\\
Z_W(\cR_W\widehat{\Phi_\lam},s)&=\lam^{2m-s}Z_W(\cR_W\widehat\Phi,s),
\end{align*}
we get
\begin{align*}
Z_W(\cR_W\Phi_\lam,2)&=\lam^{-2}Z_W(\cR_W\Phi,2),\\
Z_W(\cR_W\widehat{\Phi_\lam},2)&=\lam^{2m-2}Z_W(\cR_W\widehat\Phi,2).
\end{align*}
Then the theorem follows by integrating $\lam^sI^0(\Phi_\lam)$
over $s\in (0,1]$.
\end{proof}
Theorem \ref{rmp} in the introduction
immediately follows from the above theorem.

\bibliographystyle{plain}

\end{document}